\documentclass{article}
\usepackage{amsmath}
\newtheorem{theorem}{Theorem}
\newcommand{\mbf}[1]{\mbox{\boldmath{${#1}$}\unboldmath}}
\begin{document}
\title{Smoothed nonparametric tests and their properties}
\author{Yoshihiko MAESONO\\
Faculty of Mathematics\\
Taku MORIYAMA\\
Graduate School of Mathematics\\
Mengxin LU\\
Graduate School of Mathematics\vspace{5mm}\\
Kyushu University, Motooka 744, Fukuoka 819-0395, Japan}
\date{}

\maketitle

\begin{abstract}
In this paper we propose new smoothed sign and Wilcoxon's signed rank tests, which are based on a kernel estimator of the underlying distribution function of data.  We discuss approximations of $p$-values and asymptotic properties of these tests.  The new smoothed tests are equivalent to the ordinary sign and Wilcoxon's tests in the sense of the Pitman's asymptotic relative efficiency, and the differences of the ordinary and the new tests converge to zero in probability.  Under the null hypothesis, the main terms of the asymptotic expectations and variances of the tests do not depend on the underlying distribution.  Though the smoothed tests are not distribution-free, we can obtain Edgeworth expansions with residual term $o(n^{-1})$, which do not depend on the underlying distribution.
\end{abstract}

\section{Introduction}
\label{intro}

Let $X_1,X_2,\cdots,X_n$ be independently and identically distributed ({\it i.i.d.}) random variables with a distribution function $F(x-\theta)$, where the associated desity function satisfies $f(-x)=f(x)$, and $\theta$ is an unknown location parameter.  Here we consider a test and a confidence interval of the parameter $\theta$.  This setting is called one-sample location problem.  In order to test the null hypothesis $H_0 : \theta=0$ vs. $H_1 : \theta>0$, many nonparametric test statistics are proposed, like sign test, Wilcoxon's signed rank test etc.(see H\'{a}jek et al.~\cite{hajek1999theory})  These tests are distribution-free and have discrete distributions.  As pointed out by Lehmann \& D'Abrera~\cite{lehmann2006}, because of the discreteness of the test statistics, $p$-value jumps in response to a small change in data values, when the sample size $n$ is small or moderate.  A smoothed version of rank methods is discussed by Brown et al.~\cite{brown2001smoothed}.  They proposed a smoothed median estimator and a corresponding smoothed sign test.  The proposed test is not distribution-free and so they have discussed an Edgeworth expansion which includes unknown parameters.  Their proposed smoothed sign test has good properties, but the Pitman's asymptotic relative efficiency ($A.R.E.$) does not coincide with the ordinary sign test.  Further, when we use their Edgeworth expansion, we have to make estimators of unknown parameters.  In this paper we will propose an alternative smoothed sign test which is based on a kernel estimator of the distribution function, and discuss asymptotic properties.  We will show that its $A.R.E.$ is same as the ordinary sign test, and that the difference of two sign tests converges to zero in probability.  Further an Edgeworth expansion of the proposed sign test with residual term $o(n^{-1})$ is established and we prove the validity.  The obtained Edgeworth expansion does not depend on the underlying distribution, if we choose an appropriate kernel.

Using the same idea, we also propose a new smoothed Wilcoxon's signed rank test, and show the difference of two Wilcoxon's tests converges to zero in probability too.  An Edgeworth expansion with residual term $o(n^{-1})$ is obtained and does not depend on the underlying distribution.

Let us define $\psi(x)=1 ~(x\geq0), ~=0 ~(x<0)$ and then the sign test is equivalent to
$$S=S({\mbf X})=\sum_{i=1}^n\psi(X_i)$$
where ${\mbf X}=(X_1,X_2,\cdots,X_n)^T$.  The Wilcoxon's signed rank test is equivalent to the Mann-Whitney test
$$W=W({\mbf X})=\sum_{1\leq i\leq j\leq n}\psi(X_i+X_j).$$
For observed values ${\mbf x}=(x_1,x_2,\cdots,x_n)^T$, let us define $s=s({\mbf x})$ and $w=w({\mbf x})$ of $S$ and $W$.  If the $p$-values
$$P_0(S\geq s)\hspace{5mm}{\rm or} \hspace{5mm}P_0(W\geq w)$$
is small enough, we reject the null hypothesis $H_0$.  Here $P_0(\cdot)$ denotes a probability under the null hypothesis $H_0$.

In Table \ref{tab:table1}, based on the exact $p$-values of $S$ and $W$, we compare which $p$-value is smaller in the tail area.  Let us define
$$\Omega_{|x|}=\left\{{\mbf x}\in{\mbf R}^n ~\left\|~ |x_1|<|x_2|<\cdots<|x_n|\right.\right\}$$
and
$$\Omega_{\alpha}=\left\{{\mbf x}\in\Omega_{|x|} ~\left\|~ \frac{s-E_0(S)}{\sqrt{V_0(S)}}\geq z_{1-\alpha},\right.\hspace{4mm}{\rm or}\hspace {5mm}\frac{w-E_0(W)}{\sqrt{V_0(W)}}\geq z_{1-\alpha}\right\}$$ 
where $z_{1-\alpha}$ is a $(1-\alpha)th$ quantile of the standard normal distribution $N(0,1)$, and $E_0(\cdot)$ and $V_0(\cdot)$ are an expectation and variance under $H_0$, respectively.  The observed values $S({\mbf x})$ and $W({\mbf x})$ are invariant for the permutation of $x_1,\cdots,x_n$, and so it is sufficient to consider the case that $|x_1|<|x_2|<\cdots<|x_n|$, and $2^n$ times combinations of $sign(x_i)=\pm1 (i=1,\cdots,n)$.  We count samples that an exact $p$-value of the test is smaller than the other in the tail area $\Omega_{\alpha}$.  In Table \ref{tab:table1}, $S$ indicates a number of cases that the $p$-value of $S$ is smaller than $W$, and $W$ means a number of cases that the $p$-value of $W$ is smaller.  $W/S$ is the ratio of $W$ and $S$.  For each sample, there is one tie of the $p$-values.\\
\begin{table}
\caption{Comparison of significance probabilities \label{tab:table1}}
\begin{center}
\begin{tabular}{c|c|ccc}\hline
&sample size&$n=10$&$n=20$&$n=30$\\\hline
$z_{0.9}$&$S$&25&69080&59092679\\
&$W$&82&94442&87288529\\
&$W/S$&3.28&1.367&1.477\\\hline
$z_{0.95}$&$S$&25&32705&30857108\\
&$W$&48&47387&43957510\\
&$W/S$&1.92&1.449&1.425\\\hline
$z_{0.975}$&$S$&5&12704&14028374\\
&$W$&21&21267&22049240\\
&$W/S$&4.2&1.674&1.572\\\hline
\end{tabular}
\end{center}
\end{table}

\noindent
{\bf Remark 1}  Table \ref{tab:table1} shows that if one want to get a small $p$-value, $W$ is preferable and if not, $S$ is preferable.  This problem comes from the discreteness of the distributions of the test statistics.  In order to conquer this problem, \cite{brown2001smoothed} proposed the smoothed sign test, but its Pitman's $A.R.E.$ does not coincide with the ordinary sign test $S$.  They also obtained the Edgeworth expansion which includes unknown moments of the test statistic, and then when applying the expansion, we need estimators of them.\\

On the other hand, it is possible to use a estimator of the distribution function $F(x_0)$ as a test statistic.  Let us define the empirical distribution function
$$F_n(x_0)=\frac1n\sum_{i=1}^nI(X_i\leq x_0)$$
where $I(\cdot)$ is an indicator function.  If $F(\cdot)$ satisfies the symmetric condition, we have $F(0)=\frac12$ and so we can test the null hypothesis $H_0$ using a $p$-value
$$P_0\left\{F_n(0)\leq f_n\right\}$$
where $f_n$ is an observed value of $F_n(0)$.  This test is equivalent to the sign test $S$, that is 
\begin{equation}
S=n-nF_n(0-). \label{empiricalsign}
\end{equation}
The distribution of $F_n(0)$ is discrete and does not depend on the underlying distribution $F(\cdot)$ under $H_0$.  

In order to get a smooth estimator of the distribution function, a kernel estimator $\widetilde{F}_n$ of $F$ has been proposed.  It is natural to use $\widetilde{F}_n$ for the smoothed sign test.  Let $k(u)$ be a  kernel function which satisfies
$$\int_{-\infty}^{\infty}k(u)du=1,$$
and $K(t)$ is an integral of $k(u)$
$$K(t)=\int_{-\infty}^{t}k(u)du.$$
The kernel estimator of $F(x_0)$ is given by
$$\widetilde{F}_n(x_0)=\frac1n\sum_{i=1}^nK\Bigl(\frac{x_0-X_i}{h_n}\Bigr)$$
where $h_{n}$ is a bandwidth and $h_{n}\to0 ~(n\to\infty)$.  Thus similarly as the equation (\ref{empiricalsign}), we can use
\begin{equation}
\widetilde{S}=n-n\widetilde{F}_n(0)=n-\sum_{i=1}^nK\Bigl(-\frac{X_i}{h_n}\Bigr) \label{kernelsign}
\end{equation}
for testing $H_0$, and $\widetilde{S}$ is regarded as a smoothed version of the sign test $S$.  

The sign test $S$ or $F_n(0-)$ is distribution-free, but $\widetilde{S}$ is not.  However, under $H_0$, main terms of the asymptotic expectation and variance of $\widetilde{S}$ do not depend on $F$, i.e., asymptotically distribution-free, and we can obtain an Edgeworth expansion of $\widetilde{S}$ with residual term $o(n^{-1})$, which does not include any unknown parameters.  We will also discuss the kernel function $k(u)$ which ensures that the Edgeworth expansion does not depend on the underlying distribution $F$.\\

Similarly, we can construct the smoothed Wilcoxon's signed rank test.  Since the main term of the Mann-Whitney statistic can be regarded as an estimator of the probability
$$P\left(\frac{X_1+X_2}{2}>0\right),$$
we propose the following smoothed test statistic
$$\widetilde{W}=\frac{n(n+1)}{2}-\sum_{1\leq i\leq j\leq n}K\left(-\frac{X_i+X_j}{2h_n}\right).$$
The smoothed Wilcoxon's signed rank test $\widetilde{W}$ is not distribution-free.  However, under $H_0$, the asymptotic expectation and variance of $\widetilde{W}$ do not depend on $F$, and we can obtain an Edgeworth expansion of $\widetilde{W}$ with residual term $o(n^{-1})$.  If we use the symmetric 4-$th$ order kernel and bandwidth $h_n=o(n^{-1/4})$, the Edgeworth expansion does not depend on $F$.  Thus we can obtain the approximation of $p$-value with residual term $o(n^{-1})$.\\

In this paper we will show that the $A.R.E.$s of $\widetilde{S}$ and $\widetilde{W}$ coincide to those of $S$ and $W$, and asymptotic distributions of $S$ and $W$ are same to $\widetilde{S}$ and $\widetilde{W}$, respectively.  Further, we will show that both differences of the standardized $S$ and $\widetilde{S}$, and $W$ and $\widetilde{W}$ go to 0 in probability.  Then the smoothed test statistics are equivalent in the first order asymptotic.

In section 2, we will discuss the asymptotic properties of $\widetilde{S}$ and obtain the Edgeworth expansion with residual term $o(n^{-1})$.  A confidence interval of $\theta$ based on $\widetilde{S}$ is also discussed.  In section 3, the asymptotic properties of $\widetilde{W}$ and the Edgeworth expansion with residual term $o(n^{-1})$ will be studied, and a confidence interval of $\theta$ based on $\widetilde{W}$  is also discussed.  In section 4, we compare the obtained results by simulation.  Some proofs are given in Appendices.

\section{Asymptotic properties}
\label{sect:signtest}
In this paper we assume that the kernel $k(\cdot)$ is symmetric around the origin, that is $k(-u)=k(u)$.  Using the properties of the kernel estimator, we can obtain an expectation $E_{\theta}(\widetilde{S})$ and a variance $V_{\theta}(\widetilde{S})$.  Because of the symmetry of the underlying distribution $f$ and the kernel $k(u)$, we get
$$F(-x)=1-F(x)\hspace{5mm}{\rm and}\hspace{5mm}\int_{-\infty}^{\infty}uk(u)du=0.$$
Let us define
$$e_1(\theta)=E_{\theta}\left[1-K\left(-\frac{X_1}{h_n}\right)\right].$$
Using a transformation $u=-x/h_n$, an integration by parts and the Taylor expansion, we have
\begin{eqnarray}
e_1(\theta)&=&1-\int_{-\infty}^{\infty}K\left(-\frac{x}{h_n}\right)f(x-\theta)dx \nonumber\\
&=&1-\int_{-\infty}^{\infty} K(u)f(-\theta-h_nu)\frac{1}{h_n}du\nonumber\\
&=&1-\int_{-\infty}^{\infty} k(u)F(-\theta-h_nu)du\nonumber\\
&=&1-F(-\theta)\int_{-\infty}^{\infty} k(u)du+f(-\theta)\int_{-\infty}^{\infty}uk(u)du+O(h_n^2)\nonumber\\
&=&F(\theta)+O(h_n^2).\label{expectation}
\end{eqnarray}
Thus we get
$$E_{\theta}(\widetilde{S})=n\left\{F(\theta)+O(h_n^2)\right\}.$$
Similarly we can obtain a variance of $\widetilde{S}$.  Using the transformation $u=-x/h_n$ and Taylor expansion, we have
\begin{eqnarray*}
E_{\theta}\left[K^2\left(-\frac{X_i}{h_n}\right)\right]&=&\int_{-\infty}^{\infty} 2K(u)k(u)F(-\theta-h_nu)du\\
&=&F(-\theta)\int_{-\infty}^{\infty} 2K(u)k(u)du+O(h_n)\\
&=&F(-\theta)+O(h_n).
\end{eqnarray*}
Thus the asymptotic variance is
\begin{eqnarray*}
V_{\theta}(\widetilde{S})&=&n\left\{E_{\theta}\left[K^2\left(-\frac{X_i}{h_n}\right)\right]-[1-e_1(\theta)]^2\right\}\\
&=&n\left\{[1-F(\theta)]F(\theta)+O(h_n)\right\}.
\end{eqnarray*}

Since $\widetilde{S}$ is a sum of {\it i.i.d.} random variables, it is easy to show the asymptotic normality.  
\begin{theorem}
Let us assume that $f'$ exists and is continuous at a neighborhood of $-\theta$, and $h_n=cn^{-d} (c>0, \frac14<d<\frac12)$.  If the kernel $k(\cdot)$ is symmetric around the origin and
$$0<\lim_{n\rightarrow\infty}V_{\theta}\Bigl[1-K\Bigl(-\frac{X_i}{h_n}\Bigr)\Bigr]<\infty,$$
the standardized $\widetilde{S}$ is asymptotically normal, that is
$$\frac{\widetilde{S}-E_{\theta}(\widetilde{S})}{\sqrt{V_{\theta}(\widetilde{S})}} \hspace{5mm}\rightarrow\hspace{5mm} N(0,1)$$
in law.
\end{theorem}

Since
\begin{eqnarray*}
&&E_{\theta}(\widetilde{S})=n\left\{F(\theta)+O(h_n^2)\right\},\\
&&V_{\theta}(\widetilde{S})=n\left\{F(\theta)[1-F(\theta)]+O(h_n)\right\},\\
&&E_{\theta}(S)=nF(\theta),\\
&&V_{\theta}(S)=nF(\theta)[1-F(\theta)],
\end{eqnarray*}
we can show that the Pitman's $A.R.E.$ of $\widetilde{S}$ coincides with the sign test $S$.  Note that the main terms of the asymptotic expectation and variance of $\widetilde{S}$ do not depend on $F$ under $H_0$.

Furthermore, we can show that two sign tests are asymptotically equivalent in the sense of the first order asymptotic.  We have the following theorem.
\begin{theorem} 
Let us assume that $f'$ exists and is continuous at a neighborhood of $-\theta$, and $h_n=cn^{-d} (c>0, \frac14<d<\frac12)$.  If the kernel $k(\cdot)$ is symmetric around the origin and
$$0<\lim_{n\rightarrow\infty}V_{\theta}\Bigl[1-K\Bigl(-\frac{X_i}{h_n}\Bigr)\Bigr]<\infty,$$
the standardized sign and smoothed sign tests are equivalent, that is 
$$\lim_{n\rightarrow\infty}E_{\theta}\left\{\frac{S-E_{\theta}(S)}{\sqrt{V_{\theta}(S)}}-\frac{\widetilde{S}-E_{\theta}(\widetilde{S})}{\sqrt{V_{\theta}(\widetilde{S})}}\right\}^2=0.$$
\end{theorem}

For the sign test $S$, it is difficult to improve the normal approximation because of the discreteness of the distribution function of $S$.  The standardized sign test $S$ takes values with jump order $n^{-1/2}$, and so the formal Edgeworth expansion is meaningless.  On the other hand, since $\widetilde{S}$ is a smoothed statistic and has a continuous type distribution, we can obtain an Edgeworth expansion and prove the validity. Garc\'{i}a-Soid\'{a}n et al.~\cite{garcia1997edgeworth} discussed the Edgeworth expansion and proved the validity for the kernel estimator.  Huang \& Maesono~\cite{Huang2014improvement} have also discussed the expansion and obtained the explicit formula, for $h_n=cn^{-d} ~(c>0, ~\frac14<d<\frac12)$.  

Let us define
$$A_{i,j}=\int_{-\infty}^{\infty}K^{i}(u)k(u)u^jdu.$$
Assuming that $f'$ exists and is continuous at a neighborhood of $x_0$, and $h_n=cn^{-d} ~(c>0, ~\frac14< d<\frac12)$, Huang \& Maesono~\cite{Huang2014improvement} showed that
$$P\Biggl(\frac{\widetilde{F}_{n}(x_{0})-E[\widetilde{F}_n(x_{0})]}{\sqrt{V[\widetilde{F}_n(x_0)]}}\leq y\Biggr)=P_n(y)+o(n^{-1}),$$
where
\begin{eqnarray}
&&P_n(y)=\Phi(y)-n^{-1/2}\phi(y)Q_1(y)-n^{-1/2}h_n\phi(y)Q_1^*(y)-n^{-1}\phi(y)Q_2(y),\label{edgeworthexp}\\
&&Q_1(y)=\frac{B_{3,0}(x_0)}{6}H_{2}(y),\nonumber\\
&&Q_1^*(y)=\frac{B_{3,1}(x_0)}{6}H_{2}(y),\nonumber\\
&&Q_2(y)=\frac{B_{4,0}(x_0)}{24}H_{3}(y)-\frac{B_{3,0}^2(x_0)}{72}H_{5}(y),\nonumber\\
&&B_{3,0}(x_0)=\frac{1-2F(x_0)}{[F(x_0)\{1-F(x_0)\}]^{1/2}},\nonumber\\
&&B_{3,1}(x_0)=\frac{3f(x_0)(A_{1,1}-A_{2,1})}{[F(x_0)\{1-F(x_0)\}]^{3/2}}\nonumber\\
&&B_{4,0}(x_0)=\frac{1-3F(x_0)+3F^2(x_0)}{F(x_0)\{1-F(x_0)\}}\nonumber\\
\end{eqnarray}
and $\{H_k\}$ are the Hermite polynomials
$$H_2(y)=y^2-1,\hspace{4mm}H_3(y)=y^3-3y,\hspace{4mm}H_5(y)=y^5-10y^3-15y.$$
$\Phi(y)$ and $\phi(y)$ are distribution and density function of the standard normal $N(0,1)$. 

Since the kernel $k(u)$ is symmetric, we have
\begin{equation}
K(u)=1-K(-u),\hspace{10mm}\int_{-\infty}^{\infty}u^{2j-1}k(u)du=0 ~(j=1,2,3)\label{kernel0}
\end{equation}
and then
\begin{eqnarray*}
A_{2,1}&=&\int_{-\infty}^{\infty}K^2(u)k(u)udu=\int_{-\infty}^{\infty}\{1-K(-u)\}^2k(u)udu\\
&=&\int_{-\infty}^{\infty}k(u)udu-2\int_{-\infty}^{\infty}K(-u)k(u)udu+\int_{-\infty}^{\infty}K^2(-u)k(u)udu\\
&=&2\int_{-\infty}^{\infty}K(z)k(z)dz-\int_{-\infty}^{\infty}K^2(z)k(z)zdz=2A_{1,1}-A_{2,1}.
\end{eqnarray*}
This leads $A_{1,1}=A_{2,1}$.  Since $\widetilde{S}=n-n\widetilde{F}_n(0)$ and
$$\frac{\widetilde{S}-E(\widetilde{S})}{\sqrt{V(\widetilde{S})}}=-\frac{\widetilde{F}_n(0)-E[\widetilde{F}_n(0)]}{\sqrt{V[\widetilde{F}_n(0)]}},$$
the standardized $\widetilde{S}$ has the Edgeworth expansion with $-B_{3,0}(0), ~-B_{3,1}(0)$ and $B_{4,0}(0)$.  Further, because of $A_{1,1}=A_{2,1}$ and $F(0)=\frac12$, it is easy to see that
$$-B_{3,0}(0)=-B_{3,1}(0)=0,\hspace{4mm}B_{4,0}(0)=1.$$
Thus using Garc\'{i}a-Soid\'{a}n et al.~\cite{garcia1997edgeworth} and Huang \& Maesono~\cite{Huang2014improvement}, we have the following theorem.
\begin{theorem}
Let us assume that $k(u)$ is a symmetric kernel.  If $|f'(x)|\leq M$, $\int |u^4k(u)|du<\infty$ and the bandwidth satisfies $h_n=cn^{-d} ~(c>0, ~\frac14< d<\frac12)$, we have
$$P_0\left(\frac{\widetilde{S}-E_0(\widetilde{S})}{\sqrt{V_0(\tilde{S})}}\leq y\right)=\widetilde{P}_n(y)+o(n^{-1}),$$
where
\begin{equation}
\widetilde{P}_n(y)=\Phi(y)-\frac{1}{24n}\phi(y)H_3(y).\label{edge}
\end{equation}
\end{theorem}
The Edgeworth expansion (\ref{edge}) does not depend on the underlying distribution $F$.  In order to use this expansion we have to obtain approximations of $E_0(\widetilde{S}_n)$ and $V_0(\tilde{S})$.
\begin{theorem} 
For the positive bandwidth $h_n$, we have the following approximations of the expectation and variance under the null hypothesis $H_0$.  Here we assume that the kernel is symmetric, and $M$, $M_1$, $M_2$ and $M_3$ are some positive constants.\\
{\rm (i)} If $f^{(5)}(x)$ exists and $|f^{(5)}(x)|\leq M$ for $h_n=o(n^{-1/4})$ we have
\begin{equation}
E_0(\widetilde{S})=\frac{n}{2}+o(n^{-1/2})\label{expect}
\end{equation}
and
$$V_0(\widetilde{S})=\frac{n}{4}-2nh_nf(0)A_{1,1}-\frac{nh_n^3}{3}f''(0)A_{1,3}+o(1).$$
{\rm (ii)} If $f^{(5)}(x)$ exists and $|f^{(5)}(x)|\leq M_1$, and $A_{1,1}=A_{1,3}=0$, for $h_n=o(n^{-1/4})$ we have
\begin{equation}
V_0(\widetilde{S})=\frac{n}{4}+o(1).\label{variance}
\end{equation}
{\rm(iii)} If $f^{(4)}(x)$ exist and $|f^{(4)}(x)|\leq M_2$, and $A_{1,1}=A_{1,3}=0$, for $h_n=o(n^{-3/10})$ we have the equations {\rm(\ref{expect})} and {\rm(\ref{variance})}.\\
{\rm(iv)} If $f^{(3)}(x)$ exists and $|f^{(3)}(x)|\leq M_3$, and $A_{1,1}=0$, for $h_n=o(n^{-1/3})$ we have the equations  {\rm(\ref{expect})} and {\rm(\ref{variance})}.
\end{theorem}

\noindent
{\bf Remark 2}  In order to get above approximations, we use the Taylor expansion in the integral.  We can divide the integral at discrete points, and so we do not need to worry about the differentiability of the density function at finite points.\\

If the equations (\ref{expect}) and (\ref{variance}) hold, we can use the Edgeworth expansion of $\widetilde{S}$ for testing $H_0$ and constructing a confidence interval, without using any estimators.  For the observed value $\widetilde{s}$, we can obtain the higher order approximation of $p$-value
\begin{eqnarray}
&&P_0\left(\frac{\widetilde{S}-E_0(\widetilde{S})}{\sqrt{V_0(\widetilde{S})}}\geq \frac{\widetilde{s}-E_0(\widetilde{S})}{\sqrt{V_0(\widetilde{S})}}\right)\nonumber\\
&=&1-P_0\left(\frac{\widetilde{S}-E_0(\widetilde{S})}{\sqrt{V_0(\widetilde{S})}}\leq \frac{2}{\sqrt{n}}\left(\widetilde{s}-\frac{n}{2}\right)\right)+o(n^{-1})\nonumber\\
&=&1-\widetilde{P}_n\left(\frac{2}{\sqrt{n}}\left(\widetilde{s}-\frac{n}{2}\right)\right)+o(n^{-1}).\label{p-value}
\end{eqnarray}
Using Cornish-Fisher expansion, we can get an approximation of the $\alpha$-quantile.  Putting
$$P_0\left(\frac{\widetilde{S}-E_0(\widetilde{S})}{\sqrt{V_0(\widetilde{S})}}\leq c_{\alpha}\right)=\alpha,$$
it follows from the equation (\ref{p-value}) that
$$c_{\alpha}=z_{\alpha}+\frac{1}{24n}(z_{\alpha}^3-3z_{\alpha})+o(n^{-1})$$
where $z_{\alpha}$ is $\alpha$-quantile of $N(0,1)$.  Let us define
$$\widetilde{s}_{\alpha}=\frac{n}{2}+\frac{\sqrt{n}}{2}z_{\alpha}+\frac{1}{48\sqrt{n}}(z_{\alpha}^3-3z_{\alpha}),$$
and then we have
$$P_0(\widetilde{S}\leq\widetilde{s}_{\alpha})=\alpha+o(n^{-1}).$$
For the significance level $0<\alpha<1$, if the observed value $\widetilde{s}$ satisfies
$$\widetilde{s}\geq\widetilde{s}_{1-\alpha},$$
we reject the null hypothesis $H_0$.\\

Since the distribution of $X_i-\theta$ is the same of the distribution of $X_i$ under $H_0$, we can construct a confidence interval of $\theta$.  For the observed values ${\mbf x}=(x_1,\cdots,x_n)$, let us define
\begin{eqnarray*}
&&\widetilde{s}(\theta|{\mbf x})=n-\sum_{i=1}^nK\left(\frac{\theta-x_i}{h_n}\right),\\
&&\widehat{\theta}_U=\arg\min_{\theta}\left\{\widetilde{s}(\theta|{\mbf x})\leq\widetilde{s}_{\alpha/2}\right\}
\end{eqnarray*}
and
$$\widehat{\theta}_L=\arg\max_{\theta}\left\{\widetilde{s}_{1-\alpha/2}\leq\widetilde{s}(\theta|{\mbf x})\right\},$$
where $0<\alpha<1$.  Then the $1-\alpha$ confidence interval is given by
$$\widehat{\theta}_L\leq\theta\leq\widehat{\theta}_U.$$

\noindent
{\bf Remark 3}  The conditions of $A_{1,1}=0$ and $A_{1,3}=0$ seem restrictive, but we can construct the desired kernel.  Let us define
\begin{equation}
k(u)=\left(\frac{1}{4}(\sqrt{105}-3)+\frac12(5-\sqrt{105})|u|\right)I(|u|\leq1).\label{kernel}
\end{equation}
Then we have $A_{1,1}=0$.  $k(u)$ may take negative value, and so $\widetilde{F}_n(0)$ is not monotone increasing.  However our main purpose is to test the null hypothesis $H_0$ and to construct the confidence interval, and then we do not need to worry about it.  When we use the bandwidth $h_n=o(n^{-1/3})$, we only need the condition $A_{1,1}=0$. It is possible to construct the polynomial type kernel $k(u)$ which satisfies $A_{1,1}=A_{1,3}=0$, but it is too complicate.  Here we only consider the polynomial kernel, but it may be possible to construct another type kernel which satisfies $A_{1,1}=0$ or $A_{1,3}=0$.  We postpone this to a future work.

\section{Smoothed Wilcoxon's signed rank test}
\label{sect:wilcoxon}
Similarly as $\widetilde{S}$, we can obtain an expectation $E_{\theta}(\widetilde{W})$ and a variance $V_{\theta}(\widetilde{W})$ under $H_1$.  Let us define $g(z-\theta)$ be a density function of $\frac{X_1+X_2}{2}$.  Then we have
$$g(z-\theta)=2\int_{-\infty}^{\infty}f(u)f(2\{z-\theta\}-u)du$$
and $g(-z)=g(z)$.  Therefore, similarly as the equation (\ref{expectation}), we get the expectation
$$e_2(\theta)=E_{\theta}\left[1-K\left(-\frac{X_1+X_2}{2h_n}\right)\right]=G(\theta)+O(h_n^2),$$
where $G(z)$ is a distribution function of $\frac{X_1+X_2}{2}$ and
$$G(\theta)=\int_{-\infty}^{\infty}F(2\theta-u)f(u)du=\int_{-\infty}^{\infty}F(2\theta+u)f(u)du.$$
In order to discuss the asymptotic properties of $\widetilde{W}$, we obtain a Hoeffding\cite{hoeffding1961strong} decomposition for $U$-statistic and a representation of an asymptotic $U$-statistic, which is discussed by \cite{lai1993edgeworth}.  Let us define
\begin{eqnarray*}
\alpha_n(x)&=&2E_{\theta}\left[1-K\left(-\frac{x+X_2}{2h_n}\right)\right]-2e_2(\theta),\\
\beta_n(x,y)&=&2\left\{1-K\left(-\frac{x+y}{2h_n}\right)-e_2(\theta)-\alpha_n(x)-\alpha_n(y)\right\},\\
\alpha'_n(x)&=&2\left\{1-K\left(-\frac{x}{h_n} \right) - e_1(\theta)\right\}-2\alpha_n(x) \\
R_n&=&\frac{4}{n^{3/2}(n+1)}\sum_{i=1}^n {\alpha_n}(X_i)-\frac{2}{n^{3/2}(n+1)} \sum_{i=1}^n \alpha'_n(X_i)\\
&&-\frac{2}{n^{3/2}(n+1)} \sum_{1\leq i<j\leq n} \beta_n(X_i,X_j).
\end{eqnarray*}
Then we have
\begin{eqnarray*}
&&\frac{2}{\sqrt{n}(n+1)}\left\{\widetilde{W}-E(\widetilde{W})\right\}\\
&=&\sum_{i=1}^n \left\{ \frac{\alpha_n(X_i)}{\sqrt{n}}+\frac{\alpha'_n(X_i)}{n^{3/2}}\right\}+\sum_{1\leq i<j\leq n}\frac{\beta_n(X_i,X_j)}{n^{3/2}}+R_n.
\end{eqnarray*}
The main term of the asymptotic variance of $\widetilde{W}$ is given by $E_{\theta}[\alpha_n^2(X_1)]$.  If $|f'(x)|\leq M$, using the transformation $u=-\frac{x+y}{2h_n}$ and an integration by parts, we have
\begin{eqnarray}
&&E_{\theta}\left[1-K\left(-\frac{X_1+X_2}{2h_n}\right)\left| ~X_1=x\right.\right]\nonumber\\
&=&1-\int_{-\infty}^{\infty}K\left(-\frac{x+y}{2h_n}\right)f(y-\theta)dy\nonumber\\
&=&1-2h_n\int_{-\infty}^{\infty}K(u)f(-x-2h_nu-\theta)du\nonumber\\
&=&1-\int_{-\infty}^{\infty}k(u)F(-x-\theta-2h_nu)du\nonumber\\
&=&1-F(-x-\theta)+O(h_n^2)\nonumber\\
&=&F(x+\theta)+O(h_n^2).\label{mainapprox}
\end{eqnarray}
Here we use the fact that $\int uk(u)du=0$.  Then, if the kernel $k(\cdot)$ is symmetric, the main term of the asymptotic variance is 
\begin{eqnarray*}
E_{\theta}[\alpha_n^2(X_1)]&=&4\int_{-\infty}^{\infty}F^2(x+\theta)f(x-\theta)dx-4e_2^2(\theta)+O(h_n^2)\\
&=&4\int_{-\infty}^{\infty}F^2(u+2\theta)f(u)du-4G^2(\theta)+O(h_n^2).
\end{eqnarray*}
Using the asymptotic theory for $U$-statistics, we have the following theorem.
\begin{theorem}
Let us assume that $f'$ exists and is continuous at a neighborhood of $-\theta$, and $h_n=cn^{-d} (c>0, \frac14<d<\frac12)$.  If the kernel $k(\cdot)$ is symmetric around the origin, we have
$$\frac{\widetilde{W}-E_{\theta}(\widetilde{W})}{\sqrt{V_{\theta}(\widetilde{W})}} ~\rightarrow ~N(0,1)$$
in law.
\end{theorem}

Since
$$E_{\theta}\left[\psi(X_1+X_2)|X_1=x\right]=F(x+\theta),$$
we can show that
$$E_{\theta}\left[\frac{\{W-E_{\theta}(W)\}\{\widetilde{W}-E_{\theta}(\widetilde{W})\}}{\sqrt{V_{\theta}(W)}\sqrt{V_{\theta}(\widetilde{W})}}\right]=1+O(h_n^2)$$
and
$$\frac{V_{\theta}\left(\widetilde{W}\right)}{V_{\theta}(W)}=1+O(h_n^2).$$
Thus we have the following theorem.
\begin{theorem}
Let us assume that $f'$ exists and is continuous at a neighborhood of $-\theta$, and $h_n=cn^{-d} (c>0, \frac14<d<\frac12)$.  If the kernel $k(\cdot)$ is symmetric around the origin, the standardized $W$ and $\widetilde{W}$ are equivalent, that is 
$$\lim_{n\rightarrow\infty}E_{\theta}\left\{\frac{W-E_{\theta}(W)}{\sqrt{V_{\theta}(W)}}-\frac{\widetilde{W}-E_{\theta}(\widetilde{W})}{\sqrt{V_{\theta}(\widetilde{W})}}\right\}^2=0.$$
\end{theorem}

Bickel et al.~\cite{bickel1986edgeworth} proved the validity of the Edgeworth expansion of the $U$-statistic with residual term $o(n^{-1})$.  Since the standardized $W$ and $\widetilde{W}$ are asymptotically equivalent, modifying the results of Bickel et al.~\cite{bickel1986edgeworth}, we can obtain the Edgeworth expansion of $\widetilde{W}$ with residual term $o(n^{-1})$ and show its validity.  Similarly as the example in Bickel et al.~\cite{bickel1986edgeworth}, we can show that $\beta_n(x,y)$ satisfies the assumption in their main theorem.
\begin{theorem}
Let us assume that $f'$ exists and is continuous at a neighborhood of $-\theta$, and $h_n=cn^{-d} (c>0, \frac14<d<\frac12)$.  Then we have
\begin{equation}
P\left( \frac{2\left\{\widetilde{W}-E(\widetilde{W})\right\}}{\sqrt{n}(n+1)\xi}\leq x \right) = Q_n(x) + o(n^{-1}),\label{edgeworth2}
\end{equation}
where
\begin{eqnarray*}
Q_n(x)&=&\Phi(x)-\phi(x)\left\{\frac{P_1(x)}{\sqrt{n}}+\frac{P_2(x)}{n} \right\},\\
P_1(x)&=&\kappa_3 \frac{x^2 -1}{6 \xi^3},\\
P_2(x)&=& \left \{b_4 + \frac{b_5}{4} \right \} \frac{x}{\xi^2} + \frac{\kappa_4}{24 \xi^4} (x^3 - 3x) + \frac{\kappa_3^2 }{72\xi^6}(x^5 - 10x^3 + 15x),\\
\xi^2&=& E[\alpha_n^2(X_1)],\\
a_1&=&E[\alpha_n^3(X_1)],\\
a_2&=&E[\alpha_n(X_1) \alpha_n(X_2) \beta_n(X_1,X_2)],\\
\kappa_3&=&a_1+3a_2,\\
b_1&=&E[\alpha_n^4(X_1)],\\
b_2&=&E[\alpha_n^2(X_1)\alpha_n(X_2)\beta_n(X_1,X_2)],\\
b_3&=&E[\alpha_n(X_1)\alpha_n(X_2)\beta_n(X_1,X_3)\beta_n(X_2,X_3)],\\
b_4&=&E[\alpha_n(X_1)\alpha'_n(X_1)],\\
b_5&=&E[\beta_n^2(X_1,X_2)],\\
\kappa_4&=&b_1-3\xi^4+12(b_2+b_3).
\end{eqnarray*}
\end{theorem}

In order to use the above expansion, we have to obtain the approximation of $E_0(\widetilde{W})$ and $V_0(\widetilde{W})$.  Using the results for $U$-statistics, we have the following theorem.
\begin{theorem}
Here we assume that the kernel is symmetric.  Let $M_1$, $M_2$ and $M_3$ be some positive constants.\\
{\rm (i)} ~If $|f^{(5)}(x)|\leq M_1$ and $h_n=o(n^{-1/4})$, we have
\begin{equation}
E_0(\widetilde{W})=\frac{n(n+1)}{4}+o(n^{1/2})\label{expectwilcx}
\end{equation}
and
\begin{equation}
V_0(\widetilde{W})=\frac{n^2(2n+3)}{24}-4n^3h_n^2A_{0,2}\int_{-\infty}^{\infty}\{f(x)\}^3dx+o(n^2).\label{variancewilcx}
\end{equation}
{\rm(ii)} ~If $|f^{(4)}(x)|\leq M_2$ and  $h_n=o(n^{-3/10})$, we have the equations {\rm(\ref{expectwilcx})} and {\rm(\ref{variancewilcx})}.\\
{\rm(iii)} ~If $|f^{(3)}(x)|\leq M_3$ and $h_n=o(n^{-1/3})$ we have the equations {\rm(\ref{expectwilcx})} and {\rm(\ref{variancewilcx})}.
\end{theorem}

If the kernel is symmetric 4-$th$ order, we have $A_{0,2}=0$.  Then, under $H_0$, the Edgeworth expansion is much simplified as follows.
\begin{theorem}
Assume that $|f^{(5)}(x)|\leq M$, the kernel $k(\cdot)$ is symmetric {\rm 4}-$th$ order, and $h_n=o(n^{-1/4})$.  Then we have
\begin{equation}
P_0\left( \frac{\widetilde{W}-\frac{n(n+1)}{4}}{\sqrt{\frac{n^3}{12} +\frac{n^2}{8}}} \leq x\right)=\Phi(x)-\left(\frac{7}{20}x^3-\frac{21}{20}x\right)\label{nulapprox}
\end{equation}
\end{theorem}

The Edgeworth approximation (\ref{nulapprox}) does not depend on the underlying distribution $F$, if we use a 4-$th$ order kernel, that is
$$\int uk(u)du=\int u^2k(u)du=\int u^3k(u)du=0\hspace{3mm}{\rm and}\hspace{3mm}\int u^4k(u)du\neq0.$$

\noindent
{\bf Remark 4}  For the smoothed sign test, we have to assume $A_{1,1}=0$ when we use its Edgeworth expansion (\ref{edge}), whereas for the smoothed Wilscoxon's rank test, we only need the assumption that the kernel $k(\cdot)$ is symmetric 4-$th$ order.  Jones \& Signorini~\cite{jones1997comparison} have discussed the bias reduction method, and if we apply their method, we can easily obtain a symmetric 4-$th$ order kernel.  For the symmetric kernel $k(\cdot)$, let us define
$$k^*(u)=\frac{s_4-s_2u^2k(u)}{s_4-s_2^2}$$
where $s_2=\int u^2k(u)du$ and $s_4=\int u^4k(u)du$.  Then $k^*(u)$ is symmetric 4-$th$ order.\\

If the kernel is a symmetric 4-$th$ order, we can use the Edgeworth expansion of $\widetilde{W}$, without using any estimators.  For the observed value $\widetilde{w}$, the higher order approximation of $p$-value is given by
\begin{eqnarray}
&&P_0\left(\frac{\widetilde{W}-E_0(\widetilde{W})}{\sqrt{V_0(\widetilde{W})}}\geq \frac{\widetilde{w}-E_0(\widetilde{W})}{\sqrt{V_0(\widetilde{W})}}\right)\nonumber\\
&=&1-P_0\left(\frac{\widetilde{W}-E_0(\widetilde{W})}{\sqrt{V_0(\widetilde{W})}}\leq \frac{\widetilde{w}-\frac{n(n+1)}{4}}{\sqrt{\frac{n^3}{12} +\frac{n^2}{8}}} \right)+o(n^{-1})\nonumber\\
&=&1-\widetilde{Q}_n\left(\frac{\widetilde{w}-\frac{n(n+1)}{4}}{\sqrt{\frac{n^3}{12} +\frac{n^2}{8}}}\right)+o(n^{-1}).\label{p-valueW}
\end{eqnarray}
Using Cornish-Fisher expansion, we can get an approximation of the $\alpha$-quantile.  Putting
$$P_0\left(\frac{\widetilde{W}-E_0(\widetilde{W})}{\sqrt{V_0(\widetilde{W})}}\leq c_{\alpha}^*\right)=\alpha,$$
it follows from the equation (\ref{p-valueW}) that
$$c_{\alpha}^*=z_{\alpha}+\frac{1}{n}\left(\frac{7}{20}z_{\alpha}^3-\frac{21}{20}z_{\alpha}\right)+o(n^{-1}).$$
Let us define
$$\widetilde{w}_{\alpha}=\frac{n(n+1)}{4}+\sqrt{\frac{n^3}{12} +\frac{n^2}{8}} \left(z_{\alpha}+ \frac{1}{n}\left\{\frac{7}{20}z_{\alpha}^3-\frac{21}{20}z_{\alpha} \right\}\right),$$
we have
$$P_0(\widetilde{W}\leq\widetilde{w}_{\alpha})=\alpha+o(n^{-1}).$$
For the significance level $0<\alpha<1$, if the observed value $\widetilde{w}$ satisfies
$$\widetilde{w}\geq\widetilde{w}_{1-\alpha},$$
we reject the null hypothesis $H_0$.\\

Let us consider the confidence interval.  Since the distribution of $X_i-\theta$ is the same of the distribution of $X_i$ under $H_0$, we can use the approximation $\widetilde{Q}_n(\cdot)$.  For the observed values ${\mbf x}=(x_1,\cdots,x_n)$, let us define
$$\widetilde{w}(\theta|{\mbf x})=\frac{n(n+1)}{2}-\sum_{1\leq i\leq j\leq n}K\left(\frac{2\theta-x_i-x_j}{2h_n}\right)$$
and
\begin{eqnarray*}
&&\widehat{\theta}_U^*=\arg\min_{\theta}\left\{\widetilde{w}(\theta|{\mbf x})\leq\widetilde{w}_{\alpha/2}\right\}\\
&&\widehat{\theta}_L^*=\arg\max_{\theta}\left\{\widetilde{w}_{1-\alpha/2}\leq\widetilde{w}(\theta|{\mbf x})\right\}
\end{eqnarray*}
where $0<\alpha<1$.  Then we have the $1-\alpha$ confidence interval
$$\widehat{\theta}_L^*\leq\theta\leq\widehat{\theta}_U^*.$$

\section{Simulation study}
\label{sect:simul}

In this section, we first compare the significance probabilities of $\widetilde{S}$ and $\widetilde{W}$ by simulation.  Since the distributions of $\widetilde{S}$ and $\widetilde{W}$ depend on $F$, we compare the $p$-values by simulation.  For 100,000 times  random samples from the Normal distribution, we estimate the significance probabilities in the tail area
$$\widetilde{\Omega}_{\alpha}=\left\{{\mbf x}\in{\mbf R}^n ~\left|~ \frac{\widetilde{s}({\mbf x})-E_0(\widetilde{S})}{\sqrt{V_0(\widetilde{S})}}\geq z_{1-\alpha},\right.\hspace{4mm}{\rm or}\hspace {5mm}\frac{\widetilde{w}({\mbf x})-E_0(\widetilde{W})}{\sqrt{V_0(\widetilde{W})}}\geq z_{1-\alpha}\right\}.$$
For the simulated sample ${\mbf x}\in{\mbf R}^n$, we calculate the $p$-values based on the Edgworth expansions.  Similarly as Table \ref{tab:table1}, $\widetilde{S}$ means the significance probability of $\widetilde{S}$ is smaller than $\widetilde{W}$, and so on, in Table \ref{tab:table2}.  Comparing Table \ref{tab:table1} and \ref{tab:table2}, we can see that the differences of the $p$-values of $\widetilde{S}$ and $\widetilde{W}$ is smaller than those of $S$ and $W$.

\begin{table}
\caption{Comparison of significance probabilities \label{tab:table2}}
\begin{center}
\begin{tabular}{c|c|ccc}\hline
&sample size&$n=10$&$n=20$&$n=30$\\\hline
$z_{0.99}$&$\widetilde{S}$&7136&7716&6823\\
&$\widetilde{W}$&8174&7219&6903\\
&$\widetilde{W}/\widetilde{S}$&1.145&0.9356&1.012\\\hline
$z_{0.95}$&$\widetilde{S}$&3961&3970&3572\\
&$\widetilde{W}$&3325&3410&3331\\
&$\widetilde{W}/\widetilde{S}$&0.8394&0.8589&0.9325\\\hline
$z_{0.975}$&$\widetilde{S}$&1813&1780&1752\\
&$\widetilde{W}$&1136&1396&1555\\
&$\widetilde{W}/\widetilde{S}$&0.6266&0.7843&0.8876\\\hline
\end{tabular}
\end{center}
\end{table}

In Table \ref{tab:table3} and \ref{tab:table4}, we compare  the Edgeworth expansion by simulation.  Using the Epanechnikov kernel $k(u)=\frac{3}{4}(1-u^2)I(|u|\leq1)$ and the bandwidth $h_n=n^{-1/3}(\log n)^{-1}$, Table \ref{tab:table3} and Table \ref{tab:table4} compares simple normal approximation and the Edgeworth expansion.  Since we do not know exact distributions of the smoothed sign test $\widetilde{S}$, using the 100,000 replications of the data, we estimate values $P\left(\frac{\widetilde{S}-E_0(\widetilde{S})}{\sqrt{V_0(\widetilde{S})}}\geq z_{1-\alpha}\right)$ and denote "True" in the table.  "Edge." and "Nor." denote the Edgeworth and simple normal approximations, respectively.  The underlying distribution are the normal ($N(0,1)$), the logistic (Logis.) and the double exponential (D.Exp.).  The Epanechnikov kernel does not satisfy the condition $A_{1,1}=0$, and so we use the exact values of $f(0)$.  In Table \ref{tab:table5}, using the kernel (\ref{kernel}), which satisfies $A_{1,1}=0$, and the bandwidth $h_n=n^{-1/3}(\log n)^{-1}$, we compare the Edgeworth and normal approximation.  The double exponential distribution does not satisfy the differentiability at the origin 0, but as mentioned in Remark 1, we can show the Edgeworth expansion takes the same form.\\

\begin{table}
\caption{Edgeworth expansion ($\widetilde{S}$) for Epanechnikov kernel  \label{tab:table3}}
\begin{center}
\begin{tabular}{||c||c|c|c||c||c|c|c|c||}\hline
  $\widetilde{s}$ value& $n$=30 & $A_{1,1}\neq0$ &  & $\widetilde{s}$ value& $n$=30 &  $A_{1,1}\neq0$  &    \\
  \hline
  $z_{0.99}$ & True & Edge. & Nor. &  $z_{0.95}$ & True & Edge. & Nor. \\
  \hline
  N(0,1) & 0.00812 & \underline{0.00962} & 0.01 & N(0,1) & 0.04579 & \underline{0.04830} & 0.05 \\
  Logis. & 0.00848 & \underline{0.00962} & 0.01 & Logis. & 0.04687 & \underline{0.04829} & 0.05 \\
  D.Exp. & 0.00776 & \underline{0.00897} & 0.01 & D.Exp. & 0.04575 & \underline{0.04644} & 0.05 \\
  \hline\hline
$\widetilde{s}$ value& $n$=50 &  $A_{1,1}\neq0$  &  & $\widetilde{s}$ value& $n$=50 &  $A_{1,1}\neq0$  &   \\
  \hline
  $z_{0.99}$ & True & Edge. & Nor. & $z_{0.95}$ & True & Edge. & Nor.\\
  \hline
  N(0,1) & 0.00834 & \underline{0.00954} & 0.01 & N(0,1) & 0.04768 & \underline{0.04833} & 0.05 \\
  Logis. & 0.00836 & \underline{0.00953} & 0.01 & Logis. & 0.0481 & \underline{0.04832} & 0.05 \\
  D.Exp. & 0.00814 & \underline{0.00922} & 0.01 & D.Exp. & 0.04636 & \underline{0.04743} & 0.05 \\
  \hline\hline
$\widetilde{s}$ value& $n$=100 &  $A_{1,1}\neq0$  &  & $\widetilde{s}$ value& $n$=100 &  $A_{1,1}\neq0$  &   \\
  \hline
  $z_{0.99}$ & True & Edge. & normal & $z_{0.95}$ & True & Edge. & Nor.\\
  \hline
  N(0,1) & 0.00936 & \underline{0.00947} & 0.01 &  N(0,1) & 0.0488 & \underline{0.04835} & 0.05 \\
  Logis. & 0.00888 & \underline{0.00947} & 0.01 & Logis. & 0.04794 & \underline{0.04834} & 0.05 \\
  D.Exp. & 0.00904 & \underline{0.00946} & 0.01 & D.Exp. & 0.04803 & \underline{0.04830} & 0.05 \\
  \hline
\end{tabular}
\end{center}
\end{table}

\begin{table}
\caption{Edgeworth expansion ($\widetilde{W}$) for Epanechnikov kernel  \label{tab:table4}}
\begin{center}
\begin{tabular}{||c||c|c|c||c||c|c|c|c||}
  \hline
  $\widetilde{w}$ value  & $n$=30 & &  & $\widetilde{w}$ value  & $n$=30 & &    \\
  \hline
  $z_{0.99}$ & True & Edge. & Nor. & $z_{0.95}$ & True & Edge. & Nor.\\
  \hline
  N(0,1) & 0.00847 & 0.01123 &\underline{0.01}  & N(0,1) & 0.04891 & \underline{0.04800} & 0.05  \\
  Logis. & 0.00853 & 0.01123 &\underline{0.01}  & Logis. & 0.04635 & \underline{0.04800} & 0.05  \\
  D.Exp. & 0.00864 & 0.01123 & \underline{0.01} & D.Exp. & 0.04877 & \underline{0.04800} & 0.05  \\
  \hline
  \hline
$\widetilde{w}$ value  & $n$=50 & &  & $\widetilde{w}$ value  & $n$=50 & &   \\
  \hline
  $z_{0.99}$ & True & Edge. & Nor. & $z_{0.95}$ & True & Edge. & Nor.\\
  \hline
  N(0,1) & 0.00918 & 0.01074 &\underline{0.01}  & N(0,1) & 0.04932 & \underline{0.04880} & 0.05 \\
  Logis. & 0.00938 & 0.01074 &\underline{0.01} & Logis. & 0.04976 & 0.04880 & \underline{0.05} \\
  D.Exp. & 0.00871 & 0.01074 & \underline{0.01}  & D.exp. & 0.04874   & \underline{0.04880} & 0.05 \\
  \hline
\hline
$\widetilde{w}$ value  & $n$=100 & &  & $\widetilde{w}$ value  & $n$=100 & &   \\
  \hline
  $z_{0.99}$ & True & Edge. & Nor. & $z_{0.95}$ & True & Edge. & Nor.\\
  \hline
  N(0,1) & 0.00953  & 0.01037 & \underline{0.01} & Nor. & 0.04927 & \underline{0.04940} & 0.05 \\
  Logis. & 0.00962  & 0.01037 & \underline{0.01} & Logis. & 0.04913 & \underline{0.04940} & 0.05 \\
  D.Exp. & 0.00954 & 0.01037 & \underline{0.01} & D.Exp. & 0.04988 & 0.04940 & \underline{0.05} \\
  \hline
\end{tabular}
\end{center}
\end{table}

\begin{table}
\caption{Edgeworth expansion with the kernel $A_{1,1}=0$  \label{tab:table5}}
\begin{center}
\begin{tabular}{||c||c|c|c||c||c|c|c|c||}\hline
  $\widetilde{s}$ value& $n$=30 & $A_{1,1}=0$  &  & $\widetilde{s}$ value& $n$=30 &  $A_{1,1}=0$ &    \\
  \hline
  $z_{0.99}$ & True & Edge. & Nor. &  $z_{0.95}$ & True & Edge. & Nor. \\
  \hline
   
  N(0,1) & 0.00842 & 0.01021 &\underline{0.01} & N(0,1) & 0.05013 & 0.04993 & \underline{0.05}  \\
  Logis. & 0.00937 & 0.01021 &\underline{0.01} & Logis. & 0.0491 & \underline{0.04993} & 0.05 \\
  D.Exp. & 0.00908 & 0.01021 & \underline{0.01} & D.Exp. & 0.04903 & \underline{0.04993} & 0.05 \\
  \hline\hline
$\widetilde{s}$ value& $n$=50 & $A_{1,1}=0$ &  & $\widetilde{s}$ value& $n$=50 &  $A_{1,1}=0$ &   \\
  \hline
  $z_{0.99}$ & True & Edge. & Nor. & $z_{0.95}$ & True & Edge. & Nor.\\
  \hline
  N(0,1) & 0.0092 & 0.01012 &\underline{0.01} & N(0,1) & 0.05367 & 0.04996 & \underline{0.05} \\
  Logis. & 0.00901 & 0.01012 &\underline{0.01} & Logis. & 0.05242 & 0.04996 & \underline{0.05} \\
  D.Exp. & 0.00904 & 0.01012 & \underline{0.01} & D.Exp. & 0.05253 & 0.04996 & \underline{0.05} \\
  \hline\hline
$\widetilde{s}$ value& $n$=100 & $A_{1,1}=0$ &  & $\widetilde{s}$ value& $n$=100 &  $A_{1,1}=0$ &   \\
  \hline
  $z_{0.99}$ & True & Edge. & normal & $z_{0.95}$ & True & Edge. & Nor.\\
  \hline
  N(0,1) & 0.00962  & 0.01006 & \underline{0.01}  &  N(0,1) & 0.04903 & \underline{0.04998} & 0.05 \\
  Logis. & 0.00954  & 0.01006 & \underline{0.01} & Logis. & 0.04892 & \underline{0.04998} & 0.05 \\
  D.Exp. &0.0099 & \underline{0.01006} & 0.01 & D.Exp. & 0.04937 & \underline{0.04998} & 0.05 \\
  \hline
\end{tabular}
\end{center}
\end{table}

\noindent
{\bf Remark 5}  Azzalini~\cite{azzalini1981distribution} has recommended the bandwidth $n^{-1/3}$ for the estimation of the distribution function.  If we use the bandwidth $h_n=o(n^{-1/3})$, we only need the condition $A_{1,1}=0$.  Therefore, we will use the bandwidth $h_n=n^{-1/3}(\log n)^{-1}$ in Section \ref{sect:simul}.  Since the main term of the variance of the kernel type estimator of the distribution function does not depend on the bandwidth $h_n$, there is no trade off the bias and variance.  Then we cannot propose a clear criterion for choosing the bandwidth $h_n$.  We postpone this to a future work.\\

\noindent
{\bf Remark 6}  If we use the symmetric kernel, $n^{-1/2}$ term of the expansion is 0 and so the simple normal approximation means that the residual term is already $o(n^{-1/2})$.  The above comparisons support our results for the Edgeworth expansions.\\

Next we will compare the powers of the smoothed sign, the smoothed Wilcoxon's and the Student $t$-test.  For the significance level $\alpha$, we reject the null hypothesis $H_0$, if the observed value $\widetilde{s}$ satisfies
$$\widetilde{s}\geq \frac{n}{2}+\frac{\sqrt{n}z_{\alpha}}{2}-\frac{1}{48\sqrt{n}}(z_{1-\alpha}^3-3z_{1-\alpha}).$$
For the observed value $\widetilde{w}$, we reject $H_0$, if
$$\widetilde{w}_{\alpha}=\frac{n(n+1)}{4}+\sqrt{\frac{n^3}{12} +\frac{n^2}{8}} \left(z_{\alpha}+ \frac{1}{n}\left\{\frac{7}{20}z_{\alpha}^3-\frac{21}{20}z_{\alpha} \right\}\right).$$
The Student $t$-test statistic is given by
$$T=\frac{\sqrt{n}\overline{X}}{\sqrt{V}}$$
where
$$\overline{X}=\frac{1}{n}\sum_{i=1}^nX_i,\hspace{10mm}V=\frac{1}{n-1}\sum_{i=1}^n(X_i-\overline{X})^2.$$
For the observed value $t$, if $t\geq t(n-1;1-\alpha)$, we reject the null hypothesis $H_0$, where $t(n-1;1-\alpha)$ is a $1-\alpha$-point of $t$-distribution with $n-1$ degree of freedoms.  In Table \ref{tab:table6}, using the kernel (\ref{kernel}) and the bandwidth $h_n=n^{-1/3}(\log n)^{-1}$, we simulate the power when $\theta=0.05, ~0.1, ~0.5$ and the significance level $\alpha=0.01, ~0.05$, based on 100,000 repetitions.  In order to check the size condition, we simulate the case $\theta=0$.  When the underlying distribution $F(\cdot)$ is the double exponential, the Pitman $ARE(\widetilde{S}|T)=2, ~ARE(\widetilde{W}|T)=\frac32$, and the simulation results show that the smoothed sign test is superior than the other tests.  When the underlying distribution $F(\cdot)$ is the logistic, the Pitman $ARE(\widetilde{S}|T)=\frac{\pi^2}{12}, ~ARE(\widetilde{W}|T)=\frac{\pi^2}{9}$, and the simulation results show that the smoothed Wilcoxon's test is superior than the other tests.  The student $t$-test is superior than others, when $F(\cdot)$ is normal.  These simulation studies coincide with the comparison by the Pitman's $A.R.E.$, and so the smoothed sign and Wilcoxon's tests are exactly the continuation of the ordinal tests.  Also the simulated sizes are close to those of the significance levels.

\begin{table}
{\fontsize{10.0pt}{10.0pt}\selectfont
\caption{Power comparisons of $\widetilde{S}$, $\widetilde{W}$ and $t$-test \label{tab:table6}}
\begin{center}
\begin{tabular}{||c||c|c|c||c||c|c|c||}
   \hline
   $n$=10 & $\alpha=0.01$ & &  & $n$=10 & $\alpha=0.01$ & &  \\
  \hline
  $\theta={0}$ & $\widetilde{S}$ & $\widetilde{W}$ & $T$ & $\theta={0.05}$ & $\widetilde{S}$ & $\widetilde{W}$ & $T$ \\
  \hline
  $N(0,1)$ & 0.00908 & 0.00283 & 0.01019 & $N(0,1)$ & 0.01366 & 0.00448 & \underline{0.01435} \\
  Logis. & 0.00891  & 0.00352 & 0.00839 & Logis. & \underline{0.01369} & 0.00529  & 0.00894  \\
  D.Exp. & 0.00876 & 0.00252 & 0.00689 & D.Exp. & \underline{0.01558} & 0.00471 & 0.01102 \\
  \hline
  $\theta={0.1}$ & $\widetilde{S}$ & $\widetilde{W}$ & $T$ & $\theta={0.5}$ & $\widetilde{S}$ & $\widetilde{W}$ & $T$ \\
  \hline
  $N(0,1)$ & 0.01776 & 0.00678 & \underline{0.02034} & $N(0,1)$ & 0.12709 & 0.085 & \underline{0.16654} \\
  Logis. & \underline{0.01892} & 0.00784  & 0.01094 & Logis. & \underline{0.15677} &  0.09212 & 0.11358  \\
  D.Exp. & \underline{0.02415} & 0.00789 & 0.01795 & D.Exp. & \underline{0.23476} & 0.12665 & 0.21557 \\
  \hline
\hline
 $n$=10 & $\alpha=0.05$ & &  & $n$=10 & $\alpha=0.05$ & &     \\
  \hline
  $\theta={0}$ & $\widetilde{S}$ & $\widetilde{W}$ & $T$ & $\theta={0.05}$ & $\widetilde{S}$ & $\widetilde{W}$ & $T$ \\
  \hline
  $N(0,1)$ & 0.05151 & 0.05125 & 0.04901 & $N(0,1)$ & 0.06788 & \underline{0.06913} & 0.06804 \\
  Logis. & 0.05202 & 0.05439  & 0.0466 & Logis. & 0.06823 & \underline{0.0739}  & 0.04830  \\
  D.Exp. & 0.05153 & 0.04948 & 0.04763 & D.Exp. & \underline{0.07845} & 0.07241 & 0.06638 \\
  \hline
  $\theta={0.1}$ & $\widetilde{S}$ & $\widetilde{W}$ & $T$ & $\theta={0.5}$ & $\widetilde{S}$ & $\widetilde{W}$ & $T$ \\
  \hline
  $N(0,1)$ & 0.08451 & \underline{0.09217} & 0.08808 & $N(0,1)$ & 0.35236 & \underline{0.46168} & 0.42877 \\
  Logis. & 0.08883 & \underline{0.09846}  & 0.05615 & Logis. & 0.40752 & \underline{0.46967}  & 0.31448 \\
  D.Exp. & \underline{0.11248} & 0.10234 & 0.0945 & D.Exp. & \underline{0.52588} & 0.52544 & 0.48356 \\
  \hline\hline    
  $n$=50 & $\alpha=0.01$ & &  & $n$=50 & $\alpha=0.01$ & &   \\
  \hline
  $\theta={0}$ & $\widetilde{S}$ & $\widetilde{W}$ & $T$ & $\theta={0.05}$ & $\widetilde{S}$ & $\widetilde{W}$ & $T$ \\
  \hline
  $N(0,1)$ & 0.0088 & 0.00769 & 0.01012  & $N(0,1)$ & 0.01807 & 0.01886 & \underline{0.02385}  \\
  Logis. & 0.00829 & 0.00842 & 0.00939 & Logis. & 0.01982 & \underline{0.02112} & 0.01388\\
  D.Exp. & 0.00959 & 0.00879 & 0.00928 & D.Exp. & \underline{0.02988} & 0.02514 & 0.02331 \\
  \hline
  $\theta={0.1}$ & $\widetilde{S}$ & $\widetilde{W}$ & $T$ & $\theta={0.5}$ & $\widetilde{S}$ & $\widetilde{W}$ & $T$ \\
  \hline
  $N(0,1)$ & 0.03452 & 0.04106 & \underline{0.05034} & $N(0,1)$ & 0.65161 & 0.83087 & \underline{0.86736} \\
  Logis. & 0.04067 & \underline{0.04808} & 0.02989 & Logis. & 0.76045 & \underline{0.87510} & 0.79885\\
  D.Exp. & \underline{0.07574} & 0.06282 & 0.05212 & D.Exp. & 0.91895 & \underline{0.92874} & 0.8617 \\
  \hline
  \hline
 $n$=50 & $\alpha=0.05$ & &  & $n$=50 & $\alpha=0.05$ & &   \\
  \hline
  $\theta={0}$ & $\widetilde{S}$ & $\widetilde{W}$ & $T$ & $\theta={0.05}$ & $\widetilde{S}$ & $\widetilde{W}$ & $T$ \\
  \hline
  $N(0,1)$ & 0.05188 & 0.5068 & 0.05019 & $N(0,1)$ & 0.09057 & 0.963  & \underline{0.09749} \\
  Logis. & 0.05304  & 0.05011 & 0.04952 & Logis. & 0.09722  & \underline{0.09961} & 0.06353 \\
  D.Exp. & 0.05164 & 0.05129 & 0.05008 & D.Exp. & \underline{0.1262} & 0.1138 & 0.09917 \\
  \hline
  $\theta={0.1}$ & $\widetilde{S}$ & $\widetilde{W}$ & $T$ & $\theta={0.5}$ & $\widetilde{S}$ & $\widetilde{W}$ & $T$ \\
  \hline
  $N(0,1)$ & 0.14714 & 0.16846 & \underline{0.1719} & $N(0,1)$ & 0.87948 & 0.96125 & \underline{0.96743}  \\
  Logis. & 0.16241 & \underline{0.17977} & 0.10886 & Logis. & 0.93279 & \underline{0.974} & 0.93098 \\
  D.Exp. & \underline{0.24446} & 0.21796 & 0.18007 & D.Exp. & 0.98643 & \underline{0.98758} & 0.96177 \\
\hline
\end{tabular}
\end{center}
}
\end{table}

\section{Appendices}

{\bf Proof of Theorem 2}  For the ordinary sign test $S$, we have
$$V_{\theta}(S)=nF(\theta)[1-F(\theta)]$$
Then it is sufficient to show 
$$E_{\theta}\left[\left\{S-F(\theta)\right\}\left\{\widetilde{S}-E_{\theta}(\widetilde{S})\right\}\right]=n\left\{F(\theta)[1-F(\theta)]+O(h_n)\right\}.$$
Since $S$ and $\widetilde{S}$ are sums of $i.i.d.$ random variables, we have
\begin{eqnarray*}
&&E_{\theta}\left[\left\{S-E_{\theta}(S)\right\}\left\{\widetilde{S}-E_{\theta}(\widetilde{S})\right\}\right]\\
&=&nE_{\theta}\left[\{\psi(X_1)-E_{\theta}(\psi)\}\left\{1-K\left(-\frac{X_1}{h_n}\right)-e_1(\theta)\right\}\right].
\end{eqnarray*}
Using the transformation $u=x/h_n$, the integration by parts and the Taylor expansion, we have
\begin{eqnarray*}
&&\int_{-\infty}^{\infty}\psi(x)\left[1-K\left(-\frac{x}{h_n}\right)\right]f(x-\theta)dx\\
&=&\int_0^{\infty}\left[1-K\left(-\frac{x}{h_n}\right)\right]f(x-\theta)dx\\
&=&\int_0^{\infty}[1-K(-u)]f(h_nu-\theta)h_ndu\\
&=&\left[\{1-K(-u)\}F(h_nu-\theta)\right]_0^{\infty}-\int_0^{\infty}k(-u)F(h_nu-\theta)du\\
&=&1-\frac12F(-\theta)-F(-\theta)\int_0^{\infty}k(u)du+O(h_n)\\
&=&1-\frac12F(-\theta)-\frac12F(-\theta)+O(h_n)\\
&=&F(\theta)+O(h_n).
\end{eqnarray*}
Since $E_{\theta}(\psi)=F(\theta)$ and $E_{\theta}(1-K)=F(\theta)+O(h_n^2)$, we have
\begin{eqnarray*}
&&E_{\theta}\left[\left\{S-E_{\theta}(S)\right\}\left\{\widetilde{S}-E_{\theta}(\widetilde{S})\right\}\right]\\
&=&n\{F(\theta)-[F(\theta)]^2+O(h_n)\}.
\end{eqnarray*}
Thus we have the desired result.\\

\noindent
{\bf Proof of Theorem 4}  Assuming the differentiability of the density $f(\cdot)$, we have
\begin{eqnarray*}
\frac{1}{n}E_0(\widetilde{S})&=&1-\int_{-\infty}^{\infty} K\left(-\frac{x}{h_n}\right)f(x)dx=1-\int_{-\infty}^{\infty} k(u)F(-h_nu)du\\
&=&1-F(0)+h_nf(0)A_{0,1}-\frac{h_n^2}{2}f'(0)A_{0,2}+\frac{h_n^3}{6}f''(0)A_{0,3}\\
&&-\frac{h_n^4}{24}f^{(3)}(0)A_{0,4}+\frac{h_n^5}{120}f^{(4)}(0)A_{0,5}+O(h_n^6).
\end{eqnarray*}
Similarly, we can show that
$$E_0\left\{K^2\left(-\frac{X_1}{h_n}\right)\right\}=F(0)-2h_nf(0)A_{1,1}+h_n^2f'(0)A_{1,2}-\frac{h_n^3}{3}f''(0)A_{1,3}+O(h_n^4)$$
and
\begin{eqnarray*}
\frac{1}{n}V_0(\widetilde{S})&=&F(0)\{1-F(0)\}-2h_nf(0)A_{1,1}+h_n^2f'(0)\{A_{1,2}-F(0)A_{0,2}\}\\
&&-\frac{h_n^3}{3}f''(0)\{A_{1,3}-F(0)A_{0,3}\}+O(h_n^4).
\end{eqnarray*}
Since $k(-u)=k(u)$, we have $A_{0,1}=A_{0,3}=A_{0,5}=0$.  Further, since $f(-x)=f(x)$, we get
$$f'(x)=\lim_{\varepsilon\rightarrow0}\frac{f(x+\varepsilon)-f(x)}{\varepsilon}=-\lim_{\varepsilon\rightarrow0}\frac{f(-x-\varepsilon)-f(-x)}{-\varepsilon}=-f'(-x)$$
and then $f'(0)=0$.  Similarly, we have $f''(-x)=f''(x), ~f^{(3)}(-x)=-f^{(3)}(x)$ and $f^{(3)}(0)=0$.  Thus we have the theorem.\\

\noindent
{\bf Proof of Theorem 6}  It follows from a variance form of $U$-statistics that
$$V_{\theta}(W)=n^3\left[\int_{-\infty}^{\infty}F^2(x+2\theta)f(x)dx-G^2(\theta)\right]+O(n^2).$$
Similarly, we have
\begin{eqnarray*}
V_{\theta}\left(\widetilde{W}\right)&=&\frac{n^3}{4}E\left[\alpha_n^2(X_1)\right]+O(n^2)\\
&=&n^3\left[\int_{-\infty}^{\infty}F^2(x+2\theta)f(x)dx-G^2(\theta)\right]+O(n^2+n^3h_n^2).
\end{eqnarray*}
Thus we have
$$\lim_{n\to \infty}\frac{V_{\theta}\left(\widetilde{W}\right)}{V_{\theta}(W)}=1.$$
From direct calculations, we can show that
\begin{eqnarray*}
&&Cov_{\theta}\left(\widetilde{W}, W\right)\\
&=&n^3\left[\int_{-\infty}^{\infty}\left\{F(x+\theta)+0(h_n^2)\right\}F(x+\theta)f(x-\theta)dx-G^2(\theta)+O(h_n^2)\right]\\
&=&n^3\left[\int_{-\infty}^{\infty}F^2(x+2\theta)f(x)dx-G^2(\theta)+O(h_n^2)\right].
\end{eqnarray*}
Then we have the desired result.\\

\noindent
{\bf Proof of Theorem 8}  Here we assume that $|f^{(5)}(x)|\leq M$, and the kernel $k(\cdot)$ is symmetric.  Then we have $A_{0,1}=A_{0,3}=A_{0,5}=0$.  Since the density function $g(z)$ of $\frac{X_1+X_2}{2}$ is symmetric around the origin, similarly as $\widetilde{S}$, we get $g'(0)=g^{(3)}(0)=0$.  Therefore, we have
\begin{eqnarray*}
e_2(0)&=&E_0\left[1-K\left(-\frac{X_1+X_2}{2h_n}\right)\right]\\
&=&\frac12+h_nA_{0,1}g(0)-\frac{h_n^2}{2}A_{0,2}g'(0)+\frac{h_n^3}{6}A_{0,3}g^{(2)}(0)-\frac{h_n^4}{24}A_{0,4}g^{(3)}(0)\\
&&+\frac{h_n^5}{120}A_{0,5}g^{(4)}(0)+O(h_n^6).
\end{eqnarray*}
Similarly, as $f(x)$, we can show that $g'(0)=g^{(3)}(0)=0$, and then an approximation of the expectation of $\widetilde{W}$ under $H_0$ is given by
$$E_0[\widetilde{W}]=ne_1(0)+\frac{n(n-1)}{2}e_2(0)=\frac{n(n+1)}{4}+O(n^2h_n^6).$$
Using the transformation $u=-\frac{x+y}{2h_n}$, the integration by parts and the Taylor expansion, we get
\begin{eqnarray*}
\alpha_n(x)&=&2E_0\left[1-K\left(-\frac{X_1+X_2}{2h_n}\right) ~| ~X_1=x ~\right]-2e_1(0)\\
&=&2-2F(-x)+4h_nA_{0,1}f(-x)-2h_n^2A_{0,2}f'(-x)\\
&&+\frac{4}{3}h_n^3A_{0,3}f^{(2)}(-x)-1+O(h_n^4)\\
&=&2F(x)-1+2h_n^2A_{0,2}f'(x)+O(h_n^4).
\end{eqnarray*}
Then we have
\begin{eqnarray*}
\xi^2&=&E[\alpha_n^2(X_1)]=E_0\left[\left\{2F(X_1)-1+2h_n^2A_{0,2}f'(X_1)+O(h_n^4)\right\}^2\right]\\
&=&\frac{1}{3}+4h_n^2A_{0,2}E[\{2F(X_1)-1\}f'(X_1)]+O(h_n^4).
\end{eqnarray*}
It is easy to see that
\begin{eqnarray*}
E[F(X_1)f'(X_1)]&=&\int_{-\infty}^{\infty}F(x)f'(x)f(x)dx\\
&=&\left[\frac12F(x)\{f(x)\}^2\right]_{-\infty}^{\infty}-\frac12\int_{-\infty}^{\infty}\{f(x)\}^3dx\\
&=&-\frac12\int_{-\infty}^{\infty}\{f(x)\}^3dx
\end{eqnarray*}
and
$$E[f'(X_1)]=\int_{-\infty}^{\infty}f'(x)f(x)dx=\left[\{f(x)\}^2\right]_{-\infty}^{\infty}=0.$$
Thus we have
$$\xi^2=\frac{1}{12}-4h_n^2A_{0,2}\int_{-\infty}^{\infty}\{f(x)\}^3dx+O(h_n^4).$$

Using the representation of the variance for $U$-statistic, we have
\begin{eqnarray*}
V_0(\widetilde{W})&=&nV_0\left\{K\left(-\frac{X_1}{h_n}\right)\right\}+\frac{n(n-1)}{2}V_0\left\{K\left(-\frac{X_1+X_2}{2h_n}\right)\right\}\\
&&+2n(n-1)Cov_{0}\left\{K\left(-\frac{X_1+X_2}{2h_n}\right), K\left(-\frac{X_1}{h_n}\right)\right\}\\
&&+n(n-1)(n-2)Cov_0\left\{K\left(-\frac{X_1+X_2}{2h_n}\right), K\left(-\frac{X_1+X_3}{2h_n}\right)\right\}.
\end{eqnarray*}
Since the distribution function of $(X_1+X_2)/2$ is symmetric around the orign, we have
$$V_0\left\{K\left(-\frac{X_1}{h_n}\right)\right\}=\frac{1}{4}+O(h_n) \hspace{3mm}{\rm and}\hspace{3mm}V_0\left\{K\left(-\frac{X_1+X_2}{2h_n}\right)\right\}=\frac{1}{4}+O(h_n).$$
Furthermore, we have
\begin{eqnarray*}
&&Cov_{0}\left\{K\left(-\frac{X_1+X_2}{2h_n}\right), K\left(-\frac{X_1}{h_n}\right)\right\}\\
&=&\iint K\left(-\frac{x+y}{2h_n}\right)K\left(-\frac{x}{h_n}\right)f(x)f(y)dxdy-\frac{1}{4}+O(h_n).
\end{eqnarray*}
For the first term, using the transformation $t=x+y$, integrations by parts and Taylor expansion, we can show that
\begin{eqnarray*}
&&\iint K\left(-\frac{x+y}{2h_n}\right)K\left(-\frac{x}{h_n}\right)f(x)f(y)dxdy\\
&=&\iint K\left(-\frac{x}{h_n}\right)f(x)\left\{\left[K\left(-\frac{t}{2h_n}\right)F(t-x)\right]_{-\infty}^{\infty}\right.\\
&&\hspace*{35mm}\left.+\int k\left(-\frac{t}{2h_n}\right)F(t-x)\frac{1}{2h_n}dt\right\}dx\\
&=&\int K\left(-\frac{x}{h_n}\right)f(x)\left\{\int k(-s)F(2h_ns-x)ds\right\}dx\\
&=&\int K\left(-\frac{x}{h_n}\right)f(x)F(-x)dx+O(h_n)\\
&=&\int K\left(\frac{x}{h_n}\right)f(x)F(x)dx+O(h_n)\\
&=&\left[K\left(\frac{x}{h_n}\right)\frac{1}{2}\left\{F(x)\right\}^2\right]_{-\infty}^{\infty}-\int k\left(\frac{x}{h_n}\right)\frac{1}{2}\left\{F(x)\right\}^2\frac{1}{h_n}dx\\
&=&\frac{1}{2}-\frac{1}{2}\int k(s)\left\{F(h_ns)\right\}^2ds\\
&=&\frac{1}{2}-\frac{1}{2}\int k(s)\left\{F(0)\right\}^2ds+O(h_n)\\
&=&\frac{3}{8}+O(h_n).
\end{eqnarray*}
Thsu we have
$$Cov_{0}\left\{K\left(-\frac{X_1+X_2}{2h_n}\right), K\left(-\frac{X_1}{h_n}\right)\right\}=\frac{1}{8}+O(h_n).$$
Finally we get
\begin{eqnarray*}
&&Cov_0\left\{K\left(-\frac{X_1+X_2}{2h_n}\right), K\left(-\frac{X_1+X_3}{2h_n}\right)\right\}\\
&=&\xi^2=\frac{1}{3}-4h_n^2A_{0,2}\int_{-\infty}^{\infty}\{f(x)\}^3dx+O(h_n^4).
\end{eqnarray*}
Combining these results, we can get the variance of $\widetilde{W}$ under $H_0$.\\

\noindent
{\bf Proof of Theorem 9}  Using the transformation and the Taylor expansion, we can get approximations of $a_1, a_2, b_1,\cdots,b_5$ under $H_0$.  It follows from the approximation ({\ref{mainapprox}) that
\begin{eqnarray*}
&&E_0[\alpha_n^2(X_1)\alpha_n(X_2)\beta_n(X_1,X_2)]\\
&=&16E_0\left[\left\{F(X_1)-\frac{1}{2}\right\}^2\left\{F(X_2)-\frac{1}{2}\right\} \left\{-K\left(-\frac{X_1+X_2}{2h_n}\right)\right.\right.\\
&&\hspace{10mm}\left.\left.+1-F(X_1)-\left(F(X_2)-\frac{1}{2}\right)\right\}\right]+O(h_n^2).
\end{eqnarray*}
Using the transformation $u=\frac{x+y}{2h_n}$, the integration by parts and the Taylor expansion, we get
\begin{eqnarray*}
&&\iint \left \{ F(x) -\frac{1}{2} \right \}^2 \left \{ F(y) -\frac{1}{2} \right \}K\left(-\frac{x+y}{2h_n} \right)f(x)f(y)dxdy\\
&=&\iint \left\{F(-x+2h_nu)-\frac{1}{2}\right\}K(-u)f(-x+2h_n u)(2h_n)du\left\{ F(x)-\frac{1}{2}\right \}^2f(x)dx\\
&=&\int \left[\left\{\frac{F^2(-x +2h_n u)}{2} -\frac{F(-x +2h_n u)}{2} \right\}K(-u)\right]_{-\infty}^{\infty}\left\{F(x)-\frac{1}{2}\right\}^2 f(x)dx\\
&&+\int\left[\int_{-\infty}^{\infty}\left\{ \frac{F^2(-x +2h_n u)}{2}-\frac{F(-x +2h_n u)}{2}\right\}k(u)du\right]\left\{ F(x)-\frac{1}{2}\right \}^2f(x)dx\\
&=&\int \left \{ F(x) -\frac{1}{2} \right \}^2 \left [ \frac{F^2(-x)}{2}-\frac{F(-x)}{2}\right ] f(x) dx + O(h_n^2) \\
&=&\int_{0}^1 \left \{\frac{1}{2}-t\right \}^2 \left [ \frac{t^2}{2} -\frac{t}{2} \right ] dt + O(h_n^2) \\
&=&-\frac{1}{240}+O(h_n^2).
\end{eqnarray*}
It is easy to show
\begin{eqnarray*}
&&E_0\left[\left\{F(X_1)-\frac{1}{2}\right\}^2\left\{F(X_2)-\frac{1}{2}\right\} \left\{1-F(X_1)\right\}\right]=0,\\
&&E_0\left[\left\{F(X_1)-\frac{1}{2}\right\}^2\left\{F(X_2)-\frac{1}{2}\right\}\left\{\frac12-F(X_2)\right\}\right]=-\left(\frac{1}{12} \right)^2.
\end{eqnarray*}
Thus we have an approximation
$$E_0[\alpha_n^2(X_1)\alpha_n(X_2)\beta_n(X_1,X_2)]+O(h_n^2)=-\frac{2}{45}+O(h_n^2).$$
Similarly approximations of the moments are given by
\begin{eqnarray*}
&&E_0[\alpha_n^3(X_1)]=O(h_n^2),\\
&&E_0[\alpha_n(X_1)\alpha_n(X_2)\beta_n(X_1,X_2)]=O(h_n^2),\\
&&E_0[\alpha_n^4(X_1)]=\frac15+O(h_n^2),\\
&&E_0[\alpha_n(X_1)\alpha_n(X_2)\beta_n(X_1,X_3)\beta_n(X_2,X_3)]=\frac{2}{15}+O(h_n^2),\\
&&E_0[\alpha_n(X_1)\alpha'_n(X_1)]=-\frac16+O(h_n^2),\\
&&E_0[\beta_n^2(X_1,X_2)]=\frac13+O(h_n^2).
\end{eqnarray*}

Combining these calculations, we can get the Edgeworth expansion (\ref{nulapprox}).


\begin{thebibliography}{99}
\bibitem{azzalini1981distribution}
Azzalini, A. (1981).  A note on the estimation of a distribution funciton and quantiles by a kernel method.  {\em Biometrika} 68, 326--328.
\bibitem{bickel1986edgeworth} 
Bickel, P., G\"{o}tze, F. \& Van Zwet, W. (1986). The edgeworth expansion for U-statistics of degree two. {\em The Annals of Statistics}, 14, 1463--1484.
\bibitem{brown2001smoothed} 
Brown, B., Hall, P. \& Young, G. (2001). The smoothed median and the bootstrap. {\em Biometrika}, 88(2), 519--534.
\bibitem{garcia1997edgeworth} 
Garc\'{i}a-Soid\'{a}n, P. H., Gonz\'{a}lez-Manteiga, W. \& Prada-S\'{a}nchez, J. (1997). Edgeworth expansions for nonparametric distribution estimation with applications. {\em Journal of statistical planning and inference} 65(2), 213--231.
\bibitem{hajek1999theory}
H\'ajek, J., \v{S}id\'ak, Z. \& Sen, P. K. (1999). {\em Theory of rank tests.}, Academic Press.
\bibitem{hoeffding1961strong} 
Hoeffding, W. (1961). The strong law of large numbers for U-statistics. {\em Institute of Statistics Mimeo Series}, 302.
\bibitem{Huang2014improvement}
Huang, Z. \& Maesono, Y. (2014). Edgeworth expansion for kernel estimators of a distribution function. {\em Bulletin of Informatics and Cybernetics}, 46, 1--10.
\bibitem{jones1997comparison}
Jones, M. \& Signorini, D. (1997). A comparison of higher-order bias kernel density estimators. {\em Journal of the American Statistical Association}, 92(439), 1063--1073.
\bibitem{lai1993edgeworth}
Lai, T. L. \& Wang, J. Q. (1993). Edgeworth expansions for symmetric statistics with applications to bootstrap methods. {\em  Statistica Sinica} 3, 517--542.
\bibitem{lehmann2006}
Lehmann, E. L. \& D'Abrera, H. J. M. (2006). {\em Nonparametrics: Statistical methods based on ranks.}, Springer.

\end{thebibliography}
\end{document}